\documentclass{article}

\usepackage{amssymb, latexsym}

\usepackage{graphicx}
  
\usepackage{amsmath}
\usepackage{color}

\newtheorem{theorem}{Theorem}[section]

\newtheorem{lemma}[theorem]{Lemma}

\newtheorem{proposition}[theorem]{Proposition}

\newtheorem{remark}[theorem]{Remark}

\begin{document}
\title{ \large \bf Hamilton differential Harnack inequality and $W$-entropy for Witten Laplacian on Riemannian manifolds}
\author{\ \ Songzi Li\thanks{Research partially supported by the China Scholarship Council and a Postdoctoral Fellowship of Beijing Normal University.} , Xiang-Dong Li
\thanks{Research supported by NSFC No. 11371351, Key Laboratory RCSDS, CAS, No. 2008DP173182, and a
Hundred Talents Project of AMSS, CAS.} }

\maketitle

\begin{minipage}{120mm}
{\bf Abstract}.
In this paper, we prove the Hamilton differential Harnack inequality for positive solutions to the heat equation of the Witten Laplacian on complete Riemannian manifolds with the  $CD(-K, m)$-condition, where $m\in [n, \infty)$ and $K\geq 0$ are two constants. 
Moreover, we introduce the $W$-entropy  
and prove the $W$-entropy formula for  the fundamental solution of the Witten Laplacian on complete Riemannian manifolds with the $CD(-K, m)$-condition and on compact manifolds equipped with $(-K, m)$-super Ricci flows.

\end{minipage}


\medskip



\section{Introduction}

Differential Harnack inequality is an important tool in the study of geometric PDEs. 
Let $M$ be an $n$ dimensional complete Riemannian manifold, $u$  a
positive solution to the heat equation
\begin{eqnarray}
\partial_t u=\Delta u.\label{Heat1}
\end{eqnarray}
In their famous paper  \cite{LY}, Li and Yau proved that if $Ric\geq
-K$, where $K\geq 0$ is a positive constant, then for all
$\alpha>1$,
\begin{eqnarray}
{|\nabla u|^2\over u^2}-\alpha {\partial_t u\over u}\leq
{n\alpha^2\over 2t}+{n\alpha^2K\over \sqrt{2}(\alpha-1)}.\label{LYK}
\end{eqnarray}
In particular, if $Ric\geq 0$, then taking $\alpha\rightarrow 1$,
the Li-Yau differential Harnack inequality holds
\begin{eqnarray}
{|\nabla u|^2\over u^2}-{\partial_t u\over u}\leq {n\over
2t}.\label{LY}
\end{eqnarray}

In \cite{H1}, Hamilton proved a dimension free differential Harnack inequality on
compact Riemannian manifolds with Ricci curvature bounded from
below. More precisely, if $M$ is a compact Riemannian manifold with 
\begin{eqnarray*}
Ric\geq -K,
\end{eqnarray*}
then, for any positive and bounded solution $u$ to the heat equation
$(\ref{Heat1})$, it holds
\begin{eqnarray}
{|\nabla u|^2\over u^2}\leq \left({1\over t}+2K\right)\log(A/u),\ \ \ \forall x\in M, t>0, \label{HHA}
\end{eqnarray}
where $A:=\sup\limits\{u(t, x): x\in M, t\geq 0\}$. Indeed, the same result holds on complete Riemannian manifolds with
Ricci curvature bounded from below. Under the same condition $Ric\geq -K$, Hamilton also
proved the following   differential Harnack inequality  for any positive solution to the heat
equation $(\ref{Heat1})$
\begin{eqnarray}
{|\nabla u|^2\over u^2}-e^{2Kt}{\partial_t u\over u}\leq {n\over
2t}e^{4Kt}. \label{LYHHar}
\end{eqnarray}
In particular, when $K=0$, the above inequality reduces to the
Li-Yau Harnack inequality $(\ref{LY})$ on complete Riemannian
manifolds with non-negative Ricci curvature. 
Moreover, Hamilton \cite{H1}  proved that, on compact Riemannian manifolds with $Ric\geq -K$,  then any positive and bounded solution of the heat equation $\partial_t u=\Delta u$ with $0<u\leq A$ satisfies 
\begin{eqnarray}
{\partial_t u\over u}+{|\nabla u|^2\over u^2} \leq {K\over 1-e^{-Kt}}\left[n+4\log(A/u)\right],\ \ \ \forall\ t\geq 0.\label{HHHH2}
\end{eqnarray}

On the other hand, Perelman \cite{P1} reformulated the Ricci flow as the gradient flow of the $\mathcal{F}$-functional, where $\mathcal{F}(g, f)=\int_M (R+|\nabla f|^2)e^{-f}dv$  is defined on the product space of Riemannian metrics and $C^\infty$-functions equipped with the standard $L^2$-Riemannian metric with the constraint that $e^{-f}dv$ does not change, where $R$ is the scalar curvature of $g$.  He then introduced the $\mathcal{W}$-entropy functional  and proved its monotonicity along the conjugate equation coupled with the Ricci flow. The $\mathcal{F}$-functional  has been used by Perelman to characterize the steady gradient Ricci solitons, and the $\mathcal{W}$-entropy has been used to characterize the shrinking gradient Ricci solitons. 
As an application of the $\mathcal{W}$-entropy formula, Perelman \cite{P1} proved the non local collapsing theorem for the Ricci flow, which plays an important r\^ole for ruling out cigars, one part of the singularity
classification for the final resolution of the Poincar\'e conjecture and geometrization conjecture. 

Since Perelman's preprint \cite{P1} was posted on Arxiv in 2002, many
people have studied the $\mathcal{W}$-like entropy for other geometric
flows on Riemannian manifolds. In
\cite{N1, N2}, Ni proved the $W$-entropy formula for the 
heat equation $\partial_t u=\Delta u$ on compact and complete Riemannian manifolds with non-negative Ricci curvature, where $\Delta$ denotes the usual Laplace-Beltrami operator on Riemannian manifolds.  In \cite{LX}, Li and Xu extended Ni's $W$-entropy formula to the 
heat equation $\partial_t u=\Delta u$ on complete Riemannian manifolds with  Ricci curvature bounded from below by a negative constant.   

From \cite{P1, N1, LX, Li12, Li16, LL15}, it has been known that there is a close connection between the differential Harnack inequality and the $W$-entropy for the heat equation on Riemannian manifolds. To see this link, let $(M, g)$ be a complete Riemannian manifold with bounded geometry condition, $u$ be a positive solution to the heat equation $\partial_t u=\Delta u$. As in \cite{N1, N2}, let 
\begin{eqnarray}
H_n(u(t))=-\int_M u\log u dv-{n\over 2}(\log(4\pi t)+1).
\end{eqnarray}
Then 
\begin{eqnarray*}
{d\over dt}H_n(u(t))
=\int_M \left[{|\nabla u|^2\over u^2}-{\partial_t u\over u}-{n\over 2t}\right]u dv.
\end{eqnarray*}
Suppose that $(M, g)$ is a complete Riemannian manifold with non-negative Ricci curvature. Then the Li-Yau Harnack inequality $(\ref{LY})$ holds. This yields 
\begin{eqnarray*}
{d\over dt}H_n(u(t))\leq 0.
\end{eqnarray*}
Let
\begin{eqnarray*}
W_n(u(t))={d\over dt}(tH_n(u(t)).
\end{eqnarray*}
In \cite{N1, N2}, Ni proved that
\begin{eqnarray*}
{d\over dt}W_n(u(t))=-2t\int_M\left[ \left|\nabla^2 \log u+{g\over 2t}\right|^2+Ric(\nabla \log u, \nabla \log u)\right] udv.
\end{eqnarray*}

We now introduce some notations and definitions to develop the main part of this paper. Let $(M, g)$
be a complete Riemannian manifold, $\phi\in C^2(M)$ and $d\mu=e^{-\phi}dv$, where $v$ is the Riemannian
volume measure on $(M, g)$.  The Witten Laplacian acting on smooth functions is defined by 
$$L = e^{\phi}{\rm div}(e^{-\phi}\nabla)=\Delta - \nabla\phi\cdot\nabla.$$
For any $u, v \in C_0^\infty(M)$, the integration by parts formula holds
\begin{eqnarray*}
\int_M \langle \nabla u, \nabla v\rangle d\mu=-\int_M L u vd\mu= - \int_M u L v d\mu.
\end{eqnarray*}
Thus, $L$ is the infinitesimal generator of the Dirichlet form
\begin{eqnarray*}
\mathcal{E}(u, v)=\int_M \langle\nabla u, \nabla v\rangle d\mu, \ \ \ \ u, v\in C_0^\infty(M).
\end{eqnarray*}
In \cite{BE}, Bakry and Emery proved that for all $u\in C_0^\infty(M)$,
\begin{eqnarray}
L|\nabla u|^2-2\langle \nabla u, \nabla L u\rangle=2|\nabla^2
u|^2+2Ric(L)(\nabla u, \nabla u), \label{BWF}
\end{eqnarray}
where
$$Ric(L)= Ric + \nabla^2\phi$$
is now called the infinite dimensional Bakry-Emery Ricci curvature associated with the Witten Laplacian $L$.  For $m\in [n, \infty)$,  the $m$-dimensional Bakry-Emery Ricci curvature associated with the Witten Laplacian $L$  is defined by 
$$
Ric_{m, n}(L) = Ric + \nabla^2\phi - {\nabla\phi\otimes\nabla\phi\over m-n}.
$$
In view of this, we have (see \cite{Li05})
\begin{eqnarray*}
L|\nabla u|^2-2\langle \nabla u, \nabla L u\rangle \geq {2|Lu|^2\over m}+2Ric_{m, n}(L)(\nabla u, \nabla u). 
\end{eqnarray*}
Here we only define $Ric_{m, n}(L)$ for $m=n$ when $\phi$ is a constant. By definition, we have
$$Ric(L)=Ric_{\infty, n}(L).$$ 
Following \cite{BE}, we say that $(M, g, \mu)$ satisfies the curvature-dimension $CD(K, m)$-condition\footnote{Here the word ``$CD$'' means ``curvature-dimension''.} for a constant $K\in \mathbb{R}$ and $m\in [n, \infty]$ if and only if
$$Ric_{m, n}(L)\geq Kg.$$

Inspired by  Perelman's introduction of the modified Ricci flow $\partial_t g=-2(Ric+\nabla^2\phi)$ in \cite{P1}, we define the $(K, m)$-(Perelman) Ricci flow and $(K, m)$-super (Perelman) Ricci flows as follows. We call a manifold $(M, g(t), \phi(t), t\in [0, T])$ equipped with a family of time dependent Riemann metrics $g(t)$ and $C^2$-potentials $\phi(t)$ a $(K, m)$-(Perelman) Ricci flow if 
\begin{eqnarray*}
{1\over 2}{\partial g\over \partial t}+Ric_{m, n}(L)=Kg,      \ \ \ \forall t\in (0, T],
\end{eqnarray*}
and we call $(M, g(t), \phi(t), t\in [0, T])$  a $(K, m)$-super  (Perelman) Ricci flow if 
\begin{eqnarray*}
{1\over 2}{\partial g\over \partial t}+Ric_{m, n}(L)\geq Kg,      \ \ \ \forall t\in (0, T].
\end{eqnarray*}
See also our previous paper  \cite{LL15} and \cite{LL14, LL17a, LL17b}. 
When $\phi$ is a constant and $m=n$, the $(K, n)$-(Perelman) Ricci flow is indeed the Hamilton $K$-Ricci flow 
\begin{eqnarray*}
{1\over 2}{\partial g\over \partial t}+Ric=Kg,      \ \ \ \forall t\in (0, T],
\end{eqnarray*}
and a $(K, n)$-super  (Perelman) Ricci flow is a Hamilton $K$-super Ricci flow 
\begin{eqnarray*}
{1\over 2}{\partial g\over \partial t}+Ric\geq Kg,      \ \ \ \forall t\in (0, T].
\end{eqnarray*}
While when $m=\infty$, the $(K, \infty)$-(Perelman) Ricci flow is indeed the following extension of the modified Ricci flow introduced by Perelman \cite{P1} (where $K=0$)
\begin{eqnarray*}
{1\over 2}{\partial g\over \partial t}+Ric(L)=Kg,      \ \ \ \forall t\in (0, T],
\end{eqnarray*}
and a  $(K, \infty)$-super (Perelman) Ricci flow reads as follows 
\begin{eqnarray*}
{1\over 2}{\partial g\over \partial t}+Ric(L)\geq Kg,      \ \ \ \forall t\in (0, T].
\end{eqnarray*}
We would like to point out that the notion of super Ricci flows has been also independently introduced by K.-T. Sturm  on time-dependent metric measure spaces \cite{St3}. See also Kopfer-Sturm \cite{St4}. 

In \cite{Li05}, the Li-Yau Harnack inequality $(\ref{LYK})$  
has been extended to positive solutions of the heat equation of the Witten Laplacian 
\begin{eqnarray}
\partial_t u=Lu \label{HLu}
\end{eqnarray}
on complete Riemannian manifolds
with the $CD(-K, m)$-condition,  i.e., the $m$-dimensional Bakry-Emery Ricci curvature associated with $L$ satisfies 
$Ric_{m, n}(L)\geq -Kg$, where  $m\in [n, \infty)$ and $K\geq 0$. In particular, on complete Riemannian manifolds with the $CD(0, m)$-condition,  the classical Li-Yau Harnack inequality $(\ref{LY})$ has been extended to positive solutions to the heat equation $(\ref{HLu})$ (see \cite{Li12}) 
\begin{eqnarray}
{|\nabla u|^2\over u^2}-{\partial_t u\over u}\leq {m\over
2t}.\label{LYm}
\end{eqnarray}
In \cite{Li16}, an improved version of the Hamilton Harnack inequality $(\ref{HHA})$  has been established for  positive and bounded solutions to the heat equation $(\ref{HLu})$ on complete Riemannian manifolds
with the  $CD(-K, \infty)$-condition, where $K\geq 0$ is a constant.

On the other hand, in \cite{Li12, Li16},  the $W$-entropy formula has been also extended to the heat equation of the Witten Laplacian $(\ref{HLu})$  
on complete Riemannian manifolds with non-negative $m$-dimensional
Bakry-Emery Ricci curvature condition.  More precisely, let $(M, g)$ be a complete Riemannian manifold with bounded geometry condition, $u$ a positive solution to the heat equation $(\ref{HLu})$ of the Witten Laplacian on $(M, g, \mu)$. 
Let 
\begin{eqnarray}
H_m(u(t))=-\int_M u\log u d\mu-{m\over 2}(\log(4\pi t)+1). \label{Hm}
\end{eqnarray}
Then 
\begin{eqnarray*}
{d\over dt}H_m(u(t))
=\int_M \left[{|\nabla u|^2\over u^2}-{\partial_t u\over u}-{m\over 2t}\right]u d\mu.
\end{eqnarray*}
Suppose that $(M, g)$ is a complete Riemannian manifold with the $CD(0, m)$-condition. Then the Li-Yau Harnack inequality $(\ref{LYm})$ yields 
\begin{eqnarray*}
{d\over dt}H_m(u(t))\leq 0.
\end{eqnarray*}
Let
\begin{eqnarray}
W_m(u(t))={d\over dt}(tH_m(u(t))). \label{Wm}
\end{eqnarray}
By \cite{Li12, Li16},  we have the $W$-entropy formula for the heat equation of the Witten Laplacian 
\begin{eqnarray}
{d\over dt}W_m(u(t))&=&-2t\int_M\left[ \left|\nabla^2 \log u+{g\over 2t}\right|^2+Ric_{m, n}(L)(\nabla \log u, \nabla \log u)\right] ud\mu \nonumber\\
& &\hskip1cm -{2t\over m-n}\int_M \left(\nabla \log u\cdot \nabla \phi-{m-n\over 2t}\right)^2 ud\mu. \label{Wm00}
\end{eqnarray}
In particular, ${d\over dt}W_m(u(t))\leq 0$ if $Ric_{m, n}(L)\geq 0$. Moreover, the above definition formulas $(\ref{Hm})$ and $(\ref{Wm})$ indicate also the close  connection between the extended  Li-Yau Harnack inequality $(\ref{LYm})$ and the $W$-entropy on complete Riemannian manifolds with the $CD(0, m)$-condition.  Moreover, a rigidity theorem for $W_m$ was also proved in \cite{Li12} on complete Riemannian manifolds with the $CD(0, m)$-condition. See also \cite{Li11}  for the $W$-entropy formula for the Fokker-Planck equation on complete Riemannian manifolds with the $CD(0, m)$-condition.

In our previous papers \cite{LL15, LL17b},   we proved the $W$-entropy formula for the heat equation of the Witten Laplacian on complete Riemannian manifolds with the $CD(-K, m)$-condition, $m\in [n, \infty)$ and $K\geq 0$.  These extend Ni and Li-Xu's results from the standard case of heat equation of the Laplace-Beltrami operator on complete Riemannian manifolds with Ricci curvature condition to the general case of the heat equation of the Witten Laplacian on  complete weighted Riemannian manifolds with suitable $m$-dimensional Bakry-Emery Ricci curvature condition. 
In \cite{LL15}, when $m\in \mathbb{N}$,  we gave a direct proof of the $W$-entropy formula
 for the Witten Laplacian by applying Ni's or Li-Xu's $W$-entropy formula  for the usual Laplacian to $M\times S^{m-n}$ equipped with a suitable warped product Riemannian metric, and
 gave a natural geometric interpretation of the $W$-entropy formula  for the heat equation of  the Witten Laplacian. 
In \cite{LL15}, we  have also proved
 the $W$-entropy formula for the heat equation of time dependent Witten Laplacian on Riemannian manifolds equipped with  $(K, m)$-super Ricci flows, where $m\in [n, \infty]$ and $K\in \mathbb{R}$.  More precisely, for $K=0$,  let $(M, g(t), \phi(t), t\in [0, T])$ be 
 a compact manifolds equipped with a family of time dependent metrics $g(t)$ and $C^2$-potentials $\phi(t)$, $t\in [0, T]$ such that $d\mu=e^{-\phi}dv$ is independent of $t$ (which is equivalent to the conjugate heat equation $\partial_t \phi={1\over 2}{\rm Tr}(\partial_t g)$),  then the $W$-entropy defined by $(\ref{Hm})$ and 
 $(\ref{Wm})$ for positive solution to the heat equation $(\ref{HLu})$ of the time dependent Witten Laplacian 
$L=\Delta_{g(t)}-\nabla_{g(t)}\phi(t)\cdot\nabla_{g(t)}$ satisfies
\begin{eqnarray}
{d\over dt}W_m(u(t))&=&-2t\int_M\left[ \left|\nabla^2 \log u+{g\over 2t}\right|^2+\left({1\over 2}{\partial g\over \partial t}+Ric_{m, n}(L) \right)(\nabla \log u, \nabla \log u)\right] ud\mu\nonumber\\
& &\hskip2cm -{2t\over m-n}\int_M \left(\nabla \log u\cdot \nabla \phi-{m-n\over 2t}\right)^2 ud\mu. \label{Wmt}
\end{eqnarray}
In particular, if $(M, g(t), \phi(t), t\in [0, T])$ is a $(0, m)$-super Ricci flow in the sense that
\begin{eqnarray*}
{1\over 2}{\partial g\over \partial t}+Ric_{m, n}(L)\geq 0,\ \ \ \ {\partial \phi\over \partial t}={1\over 2}{\rm Tr}\left({\partial g\over \partial t}\right), 
\end{eqnarray*}
then $W_m(u(t))$ is decreasing in time $t$ on $[0, T]$.  For general case $K\geq 0$ and $m\in [n, \infty]$, see \cite{LL15, LL17b}.

The purpose of this paper is threefolds. First, we extend the Hamilton differential Harnack inequalities $(\ref{LYHHar})$ and $(\ref{HHHH2})$ to positive solutions of the heat equation $(\ref{HLu})$  of the Witten Laplacian on complete weighted Riemannian manifolds with the
 $CD(-K, m)$ condition. Second, we use Hamilton's Harnack inequality to introduce a new $W$-entropy and prove a new $W$-entropy formula to 
 positive solutions of the heat equation $(\ref{HLu})$  of the Witten Laplacian on complete weighted Riemannian manifolds with the
 $CD(-K, m)$ condition. We also compare our new $W$-entropy with the $W$-entropy we introduced in \cite{LL15} on complete Riemannian manifolds with the $CD(-K, m)$-condition. 
 Finally, we extend the $W$-entropy formula to the heat equation $(\ref{HLu})$ associated with the time dependent Witten-Laplacian on compact manifolds equipped with $(K, m)$-super Ricci flows.  As mentioned above, by previous works in \cite{N1, Li12, Li16, LL15},  there exists an essential link between the $W$-entropy and the Li-Yau Harnack inequality $(\ref{LYm})$ 
 for the heat equation of the Witten Laplacian on complete Riemannian manifolds satisfying the $CD(0,
m)$-condition.  Our result indicates that, when $m\in [n, \infty)$ and $K\geq 0$,  there still exists an essentially deep connection between the $W$-entropy  and the Hamilton differential Harnack inequality $(\ref{LYHHar})$ for the heat equation  of the Witten Laplacian on complete Riemannian manifolds with the $CD(-K, m)$-condition. 

\section{Statement of main results}

Now we are in a position to state our main results.  Our first result extends the Hamilton differential Harnack inequality $(\ref{LYHHar})$ to the heat equation of the Witten Laplacian on complete Riemannian manifolds with the $CD(-K, m)$-condition. 

\begin{theorem}\label{HLYH} Let $(M, g)$ be a complete Riemannian manifold and  $\phi\in C^2(M)$. Suppose that there exist some
constants $m\in [n, \infty)$ and $K\geq 0$ such that
$Ric_{m, n}(L)\geq -K$. Let $u$ be a
positive solution to the heat equation $(\ref{HLu})$. 
Then the Hamilton differential Harnack inequality holds
\begin{eqnarray}
{|\nabla u|^2\over u^2}-e^{2Kt}{\partial_t u\over u}\leq {m\over
2t}e^{4Kt}.   \label{Hamil}
\end{eqnarray}
In particular, if $Ric_{m, n}(L)\geq 0$, then the Li-Yau differential Harnack
inequality holds
\begin{eqnarray*}
{|\nabla u|^2\over u^2}-{\partial_t u\over u}\leq {m\over
2t}.
\end{eqnarray*}
\end{theorem}

Integrating the differential Harnack inequality along the geodesic on the space time,  Theorem \ref{HLYH} implies the following Harnack inequality. 

\begin{theorem}\label{Thm2}  Under the same condition and notation as in Theorem \ref{HLYH}, for all $x, y\in M$, $0<\tau<T$, we have
\begin{eqnarray*}
{u(x, \tau)\over u(y, T)}\leq \left({T\over \tau}\right)^{m/2}\exp\left\{{1\over 4}e^{2K\tau}[1+2K(T-\tau)]{d^2(x, y)\over T-\tau}+{m\over 2}[e^{2KT}-e^{2K\tau}]\right\}.
\end{eqnarray*}

\end{theorem}

The following result extends Hamilton's estimate $(\ref{HHHH2})$ to the heat equation of the Witten Laplacian on complete Riemannian manifolds with  the $CD(-K, m)$-condition.

\begin{theorem}\label{Thm3} Let $(M, g)$ be a complete Riemannian manifold with bounded Riemannian curvature tensor, $\phi\in C^2(M)$ such that $\nabla\phi$ and $\nabla^2\phi$ are uniformly bounded on $M$. Suppose that there exist some
constants $m\in [n, \infty)$ and $K\geq 0$ such that $Ric_{m, n}(L)\geq -K$. Then for any bounded and positive solution $u$ to the heat equation $(\ref{HLu})$ with  $A=\sup\limits\{u(x, t), (x, t)\in M\times [0, T]\}<\infty$, it holds
\begin{eqnarray}
{\partial_t u\over u}+{|\nabla u|^2\over u^2}\leq {K\over 1-e^{-Kt}}\left[m+4\log(A/u)\right],\ \ \ \forall\ t\in [0, T].\label{HHHH3}
\end{eqnarray}
In particular, for $t\in [0, T]$, we have
\begin{eqnarray}
{\partial_t u\over u}\leq \left(K+{1\over t}\right)\left[m+4\log(A/u)\right].\label{HHHH4}
\end{eqnarray}
\end{theorem}

We would like to mention that,  as was pointed out in the report of an anonymous referee, the above estimates are central tools in the study of the classical heat equation as well as Ricci flow and Ricci solitons, so in principle there are similar applications waiting to follow in this area for the heat equation of the Witten Laplacian as well as $(K, m)$-Ricci flow and $(K, m)$-Ricci solitons.

As an application of Theorem \ref{HLYH} and Theorem \ref{Thm3}, 
we can derive the following bound for the time derivative of the logarithm  of the heat kernel of the Witten Laplacian on complete Riemannian manifolds with the 
$CD(-K, m)$-condition, we have
\begin{eqnarray*}
-{m\over 2t}e^{2Kt}\leq \partial_t \log p_t(x, y)\leq  \left(K+{1\over t}\right)\left[m+4\log {\sup\limits_{x\in M} p_t(x, y)\over \inf\limits_{x\in M} p_t(x, y)}\right].
\end{eqnarray*}
Using the upper bound and lower bound estimates of the heat kernel $p_t(x, y)$ on complete Riemannian manifolds with the $CD(-K, m)$-condition obtained in \cite{Li05, Li12, Li16}, 
we can derive the following estimate, which 
seems new in the literature.

\begin{theorem}\label{thm4}  Under the same condition and notation as in Theorem \ref{Thm3},  for all fixed $T>0$ and $t\in (0, T]$, we have 
\begin{eqnarray*}
-{m\over 2t}e^{2Kt}\leq \partial_t \log p_t(x, y)\leq  C_{m, n, K, T}\left(1+{1\over \sqrt{t}}+{d(x, y)\over t}\right)^2, \label{htt}
\end{eqnarray*}
where $p_t(x, y)$ denotes the heat kernel of the Witten Laplacian $L$ with respect to $\mu$ on $(M, g)$. 
\end{theorem}

\begin{remark} {\rm In \cite{Li16}, it has been proved that under the condition $(M, g)$ is a complete Riemannian manifold with bounded geometry condition (i.e., the Riemannian curvature tensor as well its $k$-th covariant derivatives are uniformly 
bounded  up to the $3$-rd order), $\phi\in C^4(M)$ such that $\nabla^k\phi$ are uniformly bounded on $M$ for $1\leq k\leq 4$, then 
\begin{eqnarray*}
|\partial_t \log p_t(x, y)|\leq  C_{m, n, K, T}\left(1+{1\over \sqrt{t}}+{d(x, y)\over t}\right)^2.
\end{eqnarray*}
While Theorem \ref{thm4} need only to assume  the Riemannian curvature tensor  ${\rm Riem}$ is uniformly bounded, $Ric_{m, n}(L)\geq -Kg$ and $\phi\in C^2(M)$ such that  $\nabla\phi$ and $\nabla^2\phi$ are  uniformly bounded on $M$.}
\end{remark}

The following result indicates the  close connection between the Hamiton differential Harnack inequality $(\ref{LYHHar})$ and the $W$-entropy for the heat equation of the Witten Laplacian on complete Riemannian manifolds with the $CD(-K, m)$-condition.
When $K=0$, it reduces to the $W$-entropy formula $(\ref{Wm00})$  for the heat equation of the Witten Laplacian on complete Riemannian manifolds with the $CD(0, m)$-condition.

\begin{theorem} \label{Th-W2} Let $(M, g)$ be a complete Riemannian manifold with the bounded geometry condition and $\phi\in C^4(M)$ such that $\nabla^k\phi$ are uniformly bounded on $M$ for $1\leq k\leq 4$. Let $u$ be the heat kernel of the Witten Laplacian $L=\Delta-\nabla\phi\cdot\nabla$. Let
\begin{eqnarray*}
H_{m, K}(u, t)=-\int_M u\log u d\mu-\Phi_{m, K}(t),
\end{eqnarray*}
where $\Phi_{m, K}\in C((0, \infty), \mathbb{R})$ satisfies
\begin{eqnarray*}
\Phi_{m, K}'(t)={m\over 2t}e^{4Kt} ,\ \ \ \forall t>0.
\end{eqnarray*}
Define the $W$-entropy by the Boltzmann formula
\begin{eqnarray*}
W_{m, K}(u, t)={d\over dt}(tH_{m, K}(u, t)).
\end{eqnarray*}
Then
\begin{eqnarray*}
{d\over dt}W_{m, K}(u, t)
&=&-2t\int_M \left|\nabla^2\log u+\left({K\over 2}+{1\over 2t}\right)g\right|^2 ud\mu\\
& &\hskip0.5cm -2t\int_M \left(Ric_{m, n}(L)+Kg\right)(\nabla\log u, \nabla \log u) ud\mu\\
& &\hskip1cm -{2t\over m-n}\int_M \left|\nabla \phi\cdot \nabla\log  u-{(m-n)(1+Kt)\over 2t}\right|^2ud\mu\\
& &\hskip2cm -{m\over 2t}\left[e^{4Kt}(1+4Kt)-(1+Kt)^2\right].
\end{eqnarray*}
In particular, if $Ric_{m, n}(L)\geq -Kg$, then, for all $t>0$, we have
\begin{eqnarray*}
{d\over dt}W_{m, K}(u, t)\leq -{m\over
2t}\left[e^{4Kt}(1+4Kt)-(1+Kt)^2\right].
\end{eqnarray*}
Moreover, the equality holds at some time $t=t_0>0$ if and only if
$(M, g, \phi)$ is a  fixed point of the $(-K, m)$-Ricci flow, called $(-K, m)$-Ricci soliton or $(-K, m)$-quasi-Einstein manifold, 
$$Ric_{m, n}(L)=-Kg,$$ the
potential function $f=-\log u$ satisfies the shrinking soliton
equation with respect to $Ric_{m, n}(L)$, i.e.,
\begin{eqnarray*}
Ric_{m, n}(L)+2\nabla^2f={g\over t},
\end{eqnarray*}
and moreover
\begin{eqnarray*}
\nabla \phi\cdot \nabla f=-{(m-n)(1+Kt)\over 2t}.
\end{eqnarray*}
\end{theorem}

The following result extends Theorem \ref{Th-W2} to the heat equation of the time dependent Witten Laplacian on compact manifolds equipped with $(-K, m)$-super Ricci flows. When $K=0$, it is the  $W$-entropy formula  $(\ref{Wmt})$ on the $(0, m)$-super Ricci flows, which was proved in our previous paper \cite{LL15}.

\begin{theorem} \label{Th-W3} Let $(M, g(t), \phi(t), t\in [0, T])$ be a compact manifold with a family of Riemannian metrics and $C^\infty$-potentials $(g(t), \phi(t), t\in [0, T])$ . Suppose that 
\begin{eqnarray*}
{\partial \phi\over \partial t}={1\over 2}{\rm Tr}\left({\partial g\over \partial t}\right).
\end{eqnarray*}
Let $u$ be a positive solution to the heat equation $(\ref{HLu})$  of the time dependent Witten Laplacian 
$$L=\Delta_{g(t)}-\nabla_{g(t)}\phi(t)\cdot\nabla_{g(t)}.$$ 
Let $H_{m, K}(u, t)$  and $W_{m, K}(u, t)$ be as in Theorem \ref{Th-W2}. 
Then
\begin{eqnarray*}
{d\over dt}W_{m, K}(u, t)
&=&-2t\int_M \left|\nabla^2\log u+\left({K\over 2}+{1\over 2t}\right)g\right|^2ud\mu\\
& &\hskip0.5cm -2t\int_M \left({1\over 2}{\partial g\over \partial t}+Ric_{m, n}(L)+Kg\right)(\nabla\log u, \nabla \log u)ud\mu\\
& &\hskip1cm -{2t\over m-n}\int_M \left|\nabla \phi\cdot \nabla\log  u-{(m-n)(1+Kt)\over 2t}\right|^2ud\mu\\
& &\hskip1.5cm -{m\over 2t}\left[e^{4Kt}(1+4Kt)-(1+Kt)^2\right].
\end{eqnarray*}
In particular, if $(M, g(t), \phi(t), t\in [0, T])$ is a compact manifolds equipped with a $(-K, m)$-super Ricci flow in the sense that
\begin{eqnarray*}
{1\over 2}{\partial g\over \partial t}+Ric_{m, n}(L)\geq -K, \hskip0.3cm {\partial \phi\over \partial t}={1\over 2}{\rm Tr}\left({\partial g\over \partial t}\right),
\end{eqnarray*}
then for all $t\in (0, T]$,
we have
\begin{eqnarray*}
{d\over dt}W_{m, K}(u, t)\leq -{m\over
2t}\left[e^{4Kt}(1+4Kt)-(1+Kt)^2\right].
\end{eqnarray*}
Moreover, the equality holds on $(0, T]$ if and only if
$(M, g(t), \phi(t), t\in [0, T])$ is a $(-K, m)$-Ricci flow  in the sense that
\begin{eqnarray*}
{\partial g\over \partial t}&=&-2(Ric_{m, n}(L)+Kg), \nonumber\\
{\partial \phi\over \partial t}&=&-R-\Delta \phi-{|\nabla \phi|^2\over m-n}-nK,
\end{eqnarray*}
 the
potential function $f=-\log u$ satisfies the Hessian equation
\begin{eqnarray*}
\nabla^2f=\left({K\over 2}+{1\over 2t}\right)g,
\end{eqnarray*}
and moreover
\begin{eqnarray*}
\nabla \phi\cdot \nabla f=-{(m-n)(1+Kt)\over 2t}.
\end{eqnarray*}
\end{theorem}
%
%
%
%

We can also extend  the Hamilton Harnack inequalities to positive solutions to the heat equation $\partial_t u=Lu$ associated with the time dependent Witten Laplacian $L=\Delta-\nabla\phi\cdot\nabla$ on compact or complete  Riemannian manifolds equipped with a variant 
of $(-K, m)$-super Ricci flow.  To save  the length of the paper, we will do it in a forthcoming paper. See \cite{LL14}.

The rest of this paper is organized as follows.  In Section $3$, we prove   Theorem \ref{HLYH}, Theorem \ref{Thm2} and Theorem \ref{Thm3}.  In Section $4$,  we prove  Theorem \ref{Th-W2} and Theorem \ref{Th-W3}. In Section $5$, 
we compare the $W$-entropy in Theorem \ref{Th-W2} and Theorem \ref{Th-W3} with the $W$-entropy defined in our previous paper \cite{LL15}.    

This paper is an improved version of a part of our previous preprint \cite{LL14}.  Due to the limit of the length of the paper, we split \cite{LL14} into several papers. See also \cite{LL17a, LL17b, LL17c}. 

\section{Hamilton Harnack inequalities for Witten Laplacian}

\subsection{Proof of Theorem \ref{HLYH}}

By the generalized Bochner-Weitzenb\"ock formula, we have
\begin{eqnarray*}
(L-\partial_t){|\nabla u|^2\over u}={2\over u}\left|\nabla^2
u-{\nabla u\otimes \nabla u\over u}\right|^2+{2\over u}Ric(L)(\nabla
u, \nabla u).\label{GB2}
\end{eqnarray*}
Taking trace in the first quantity on the right hand side, we can
derive
\begin{eqnarray*}
(L-\partial_t){|\nabla u|^2\over u}\geq {2\over nu}\left|\Delta
u-{|\nabla u|^2\over u}\right|^2+{2\over u}Ric(L)(\nabla u, \nabla
u).\label{GB3}
\end{eqnarray*}
Applying the inequality
\begin{eqnarray*}
(a+b)^2\geq {a^2\over 1+\alpha}-{b^2\over \alpha}, \ \ \forall \alpha>0, 
\end{eqnarray*}
to $a=\partial_t u-{|\nabla u|^2\over u}$, $b=\nabla \phi\cdot\nabla u$, and $\alpha={m-n\over n}$, we have
\begin{eqnarray*}
(L-\partial_t){|\nabla u|^2\over u}\geq {2\over m u}\left|\partial_t
u-{|\nabla u|^2\over u}\right|^2+{2\over u} Ric_{m, n}(L)(\nabla u,
\nabla u). \label{GB4}
\end{eqnarray*}
Hence, under the condition $Ric_{m, n}(L)\geq -K$, it holds
\begin{eqnarray*}
(L-\partial_t){|\nabla u|^2\over u} \geq  {2\over m
u}\left|\partial_t u-{|\nabla u|^2\over u}\right|^2-{2K|\nabla
u|^2\over u}.\label{GB5}
\end{eqnarray*}
Let
\begin{eqnarray*}
h={\partial u \over \partial t}-e^{-2Kt}{|\nabla u|^2\over u}+e^{2Kt}{m\over 2t}u.
\end{eqnarray*}
Then $\lim\limits_{t\rightarrow 0^+}h(t)=+\infty$, and
\begin{eqnarray}
(\partial_t-L)h\geq {2\over m u}e^{-2Kt}\left|\partial_t u-{|\nabla u|^2\over u}\right|^2-e^{2Kt}{m\over 2t^2}u. \label{ht}
\end{eqnarray}

We now prove that $h\geq 0$ on $M\times \mathbb{R}^+$. In compact case, suppose that $h$ attends its minimum at some $(x_0, t_0)$ and $h(x_0, t_0)<0$. Then, at $(x_0, t_0)$, it holds
\begin{eqnarray*}
{\partial h\over \partial t}\leq 0,\ \ \Delta h\geq 0, \ \nabla h=0.
\end{eqnarray*}
Thus at $(x_0, t_0)$, $(\partial_t-L)h\leq 0$.  On the other hand,  as $h(x_0, t_0)<0$, we have
\begin{eqnarray*}
0\leq e^{2Kt}{m\over 2t}u<e^{-2Kt}{|\nabla u|^2\over u}-{\partial u\over \partial t}\leq {|\nabla u|^2\over u}-{\partial u\over \partial t},
\end{eqnarray*}
and hence by $(\ref{ht})$ we have
\begin{eqnarray*}
(\partial_t-L) h>0.
\end{eqnarray*}
This finishes the proof of Theorem \ref{HLYH} in compact case.

In complete non-compact case, let $f=\log u$, and let
\begin{eqnarray*}
F=te^{-2Kt}(e^{-2Kt}|\nabla f|^2-f_t)=te^{-4Kt}|\nabla f|^2-te^{-2Kt}f_t.
\end{eqnarray*}
Obviously, $F(0, x)\equiv 0$. We shall prove that
\begin{eqnarray*}
F\leq {m\over 2}.
\end{eqnarray*}
By direct calculation
\begin{eqnarray*}
LF&=&te^{-4Kt}L|\nabla f|^2-te^{-2Kt}Lf_t\\
\partial_t F
&=&(1-4Kt)e^{-4Kt}|\nabla f|^2+(2Kt-1)e^{-2Kt}f_t+te^{-4Kt}\partial_t |\nabla f|^2-te^{-2Kt}f_{tt},
\end{eqnarray*}
we have
\begin{eqnarray*}
(L-\partial_t)F&=&te^{-4Kt}(L-\partial_t)|\nabla f|^2-te^{-2Kt}(L-\partial_t)f_t\\
& &\ \ \ +(4Kt-1)e^{-4Kt}|\nabla f|^2-(2Kt-1)e^{-2Kt}f_t.
\end{eqnarray*}
By the generalized Bochner formula, it holds
\begin{eqnarray*}
(L-\partial_t)|\nabla f|^2=2|\nabla^2 f|^2+2Ric(L)(\nabla f, \nabla f)-4\nabla^2f(\nabla f, \nabla f).
\end{eqnarray*}
Note that
\begin{eqnarray*}
Lf_t&=&L\left({Lu\over u}\right)
={L^2u\over u}-2\langle\nabla Lu, {\nabla u\over u^2}\rangle+Lu\left(-{Lu\over u^2}+2{|\nabla u|^2\over u^3}\right),\\
\partial_t f_t&=&\partial_t\left({Lu\over u}\right)={L^2u\over u}-{|Lu|^2\over u^2},
\end{eqnarray*}
which yields
\begin{eqnarray*}
(L-\partial_t)f_t&=& 2{Lu|\nabla u|^2\over u^3} -2\langle\nabla Lu, {\nabla u\over u^2}\rangle\\
&=&-4\nabla^2 f(\nabla f, \nabla f)-2\langle\nabla Lf, \nabla f\rangle.
\end{eqnarray*}
Hence
\begin{eqnarray*}
(L-\partial_t)F
&=& 2te^{-4Kt}[|\nabla^2 f|^2 + 2(e^{2Kt}-1)\nabla^2 f(\nabla f, \nabla f)]\\
& &+2te^{-4Kt}Ric(L)(\nabla f, \nabla f) +2 te^{-2Kt}\langle\nabla Lf, \nabla f\rangle\\
& &\ \ \ +(4Kt-1)e^{-4Kt}|\nabla f|^2-(2Kt-1)e^{-2Kt}(Lf+|\nabla f|^2).
\end{eqnarray*}
Now
\begin{eqnarray*}
F &=&  te^{-4Kt}(1 - e^{2Kt})|\nabla f|^2 - te^{-2Kt}Lf,\\
\langle \nabla F, \nabla f \rangle &=& 2te^{-4Kt}(1 - e^{2Kt})\nabla^2f(\nabla f, \nabla f) - te^{-2Kt}\langle\nabla Lf, \nabla f\rangle.
\end{eqnarray*}
Therefore
\begin{eqnarray*}
(L-\partial_t)F
&=&  2te^{-4Kt}|\nabla^2 f|^2 -2\langle \nabla F, \nabla f \rangle\\
& &\ \ +2te^{-4Kt}\left(Ric(L)(\nabla f, \nabla f)+K|\nabla
f|^2\right)+ \frac{(2Kt-1)}{t}F.
\end{eqnarray*}
Note that
$$
|\nabla^2 f|^2 \geq \frac{1}{n}|\Delta f|^2 \geq \frac{1}{m}|L f|^2
-{1\over m-n}\nabla\phi\otimes\nabla \phi(\nabla f, \nabla f).
$$
Thus
\begin{eqnarray*}
(L-\partial_t)F&\geq& 2te^{-4Kt}{|Lf|^2\over m}-2\langle \nabla F, \nabla f \rangle\\
& &\ \ +2te^{-4Kt}\left(Ric_{m, n}(L)(\nabla f, \nabla f)+K|\nabla
f|^2\right)+ \frac{(2Kt-1)}{t}F\\
&\geq& \frac{2te^{-4Kt}}{m}\left[\frac{(te^{-2Kt}(e^{-2Kt}
-1)|\nabla f|^2 - F)^2}{t^2e^{-4Kt}}\right]-2\langle \nabla F, \nabla f \rangle +  \frac{(2Kt-1)}{t}F\\
&\geq& \frac{2[te^{-2Kt}(e^{-2Kt} -1)|\nabla f|^2 - F]^2}{mt}-
2\langle \nabla F, \nabla f \rangle +  \frac{(2Kt-1)}{t}F.
\end{eqnarray*}

Similarly to \cite{Li05}, let $\eta$ be a $C^2$-function on $[0,\infty)$ such that $\eta=1$ on $[0, 1]$ and $\eta=0$ on $[2,\infty)$, with $-C_1\eta^{1/2}(r)\leq \eta'(r)\leq 0$, and $\eta''(r)\geq C_2$, where $C_1>0$ and $C_2>0$ are two constants. Let $\rho(x)=d(o, x)$ and define $\psi(x)=\eta({\rho(x)/R})$.
Since $\rho$ is Lipschitz on the complement of the cut locus of $o$, $\psi$ is a Lipschitz function with support
in $B(o, 2R)\times [0,\infty)$. As explained in Li and Yau \cite{LY}, an argument of Calabi allows us
to apply the maximum principle to $\psi F$. Let $(x_0, t_0)\in M\times [0, T]$ be
a point where $\psi F$ achieves the maximum. Then, at $(x_0, t_0)$,
\begin{eqnarray*}
\partial_t(\psi F)\geq 0, \ \Delta(\psi F)\leq 0, \ \nabla(\psi F)=0.
\end{eqnarray*}
This yields
\begin{eqnarray*}
(L-\partial_t)(\psi F)=\Delta (\psi F)-\nabla\phi\cdot\nabla(\psi
F)-\partial_t(\psi F)\leq 0.
\end{eqnarray*}
Similarly to \cite{Li05}, we have
\begin{eqnarray*}
(L-\partial_t)(\psi F)&=&\psi(L-\partial_t)F+(L\psi)F+2\nabla\psi \cdot\nabla F\\
&\geq&\psi(L-\partial_t)F-A(R)F+2\nabla\psi\cdot\nabla F\\
&\geq &\psi(L-\partial_t)F -A(R)F+2\langle\nabla \psi, \nabla(\psi
F)\rangle\psi^{-1}-2F|\nabla\psi|^2\psi^{-1}.
\end{eqnarray*}where we use
\begin{eqnarray*}
L\psi\geq -A(R):=-{C_1\over
R}(m-1)\sqrt{K}\coth(\sqrt{K}R)-{C_2\over R^2},
\end{eqnarray*}
and for some constant $C_3>0$
$${|\nabla \psi|^2\over \psi}\leq {C_3\over R^2}.$$
Let $C(n, K, R)={C_1\over
R}(m-1)\sqrt{K}\coth(\sqrt{K}R)+{C_2+C_3\over R^2}$. At the point $(x_0, t_0)$, we have
\begin{eqnarray*}
0 &\geq &\psi(L-\partial_t)F -(A(R) + 2|\nabla\psi|^2\psi^{-1})F\\
&\geq & \psi\left[\frac{2[te^{-2Kt}(e^{-2Kt} -1)|\nabla f|^2 - F]^2}{mt}- 2\langle \nabla F, \nabla f \rangle +  \frac{(2Kt-1)}{t}F\right] -C(n, K, R)F\\
&\geq& \psi\frac{2}{mt}F^2 + \psi\frac{4e^{-2Kt}(1 - e^{-2Kt})|\nabla f|^2}{m}F + 2F\langle \nabla \psi, \nabla f \rangle +  \left[(2K -\frac{1}{t})\psi -C(n, K, R)\right]F\\
&\geq& \psi\frac{2}{mt}F^2 + \psi\frac{4e^{-2Kt}(1 - e^{-2Kt})|\nabla f|^2}{m}F - 2F|\nabla \psi||\nabla f| +  \left[(2K -\frac{1}{t})\psi -C(n, K, R)\right]F\\
&\geq& \psi\frac{2}{mt}F^2 + \psi\frac{4e^{-2Kt}(1 - e^{-2Kt})|\nabla f|^2}{m}F - 2{C_2 \over R}F\psi^{1/2}|\nabla f| +  \left[(2K -\frac{1}{t})\psi -C(n, K, R)\right]F.
\end{eqnarray*}
Multiplying by $t$ on both sides, and using the Cauchy-Schwartz inequality, we get
\begin{eqnarray*}
0 
&=& \psi\frac{2}{m}F^2 + tF\left[\psi\frac{4e^{-2Kt}(1 - e^{-2Kt})|\nabla f|^2}{m} - 2{C_2 \over R}\psi^{1/2}|\nabla f|\right] +  [(2Kt - 1)\psi - C(n, K, R)t]F\\
&\geq& \psi\frac{2}{m}F^2 +  \left[(2Kt - 1)\psi - C(n, K, R)t -
{C_2 mt\over 4e^{-2Kt}(1 - e^{-2Kt})R^2}\right]F.
\end{eqnarray*}
Notice that the above calculation is done at the point $(x_0, t_0)$.
Since $\psi F$ reaches its maximum at this point, we can assume that
$\psi F(x_0, t_0) > 0$. Thus
\begin{eqnarray*}
0 
\geq \frac{2}{m}(\psi F)^2  -\left[1 + C(n, K, R)t + {C_2 m\over
4e^{-2Kt}(1 - e^{-2Kt})R^2}t\right](\psi F),
\end{eqnarray*}
which yields that, for any $(x, t) \in B_{R} \times [0, T]$,
\begin{eqnarray*}
F(x,t) &\leq& (\psi F)(x_0, t_0) \leq \frac{m}{2}\left[1 + C(n, K, R)t_0 + {C_2m \over 4e^{-2Kt_0}(1 - e^{-2Kt_0})R^2}t_0\right]\\
&\leq& \frac{m}{2} \left[1 + C(n, K, R)T + \max\limits_{t\in [0,
T]}{C_2 mt\over 4e^{-2Kt}(1 - e^{-2Kt})R^2}\right].
\end{eqnarray*}
Let $R \rightarrow \infty$, we obtain
$$
F \leq \frac{m}{2}.
$$
The proof of Theorem \ref{HLYH} is completed. \hfill $\square$

\medskip

\subsection{Proof of Theorem \ref{Thm2}}.  

The proof is as the same as the one of Corollary 2.2 in \cite{H1}. For the completeness we reproduce it as follows. Let $l(x, t)=\log u(x, t)$. Then the Hamilton Harnack inequality is equivalent to
\begin{eqnarray}
{\partial l\over \partial t}-e^{-2Kt}|\nabla l|^2+e^{2Kt}{m\over 2t}\geq 0.\label{HHHH}
\end{eqnarray}
Let $\gamma: [0, T]\rightarrow M$ be a geodesic with reparametrization by arc length $s: [\tau, T]\rightarrow [0, T]$ so that $\gamma(s(\tau))=x$ and $\gamma(s(T))=y$. Let $S(t)={d\gamma(s(t))\over dt}=\dot\gamma(s(t))\dot s(t)$. Then
$|\dot\gamma(s(t))|=1$. Integrating along $\gamma(s(t))$ from $t=\tau$ to $t=T$, we have
\begin{eqnarray*}
l(y, T)-l(x, \tau)=\int_\tau^T \left[{\partial l\over \partial t}+\nabla l\cdot S\right]dt.
\end{eqnarray*}
By the Cauchy-Schwartz inequality
\begin{eqnarray*}
e^{-2Kt}|\nabla l|^2+{1\over 4}e^{2Kt}|S|^2\geq \nabla l\cdot S
\end{eqnarray*}
From this and $(\ref{HHHH})$ we obtain
\begin{eqnarray*}
l(y, T)-l(x, \tau)\geq -{1\over 4}\int_\tau^Te^{2Kt}|S|^2dt-\int_\tau^T {m\over 2t}e^{2Kt}dt.
\end{eqnarray*}
Note that $d(x, y)=\int_\tau^T |S|dt=\int_\tau^T ds(t)$.
Choosing $s(t)=a[e^{-2K\tau}-e^{-2Kt}]$, with
$$a={d(x, y)\over e^{-2K\tau}-e^{-2KT}},$$ we have
\begin{eqnarray*}
\int_\tau^Te^{2Kt}|S|^2dt=\int_t^T e^{2Kt}{\dot s}^2(t)dt={2Kd^2(x, y)\over e^{-2K\tau}-e^{-2KT}}.
\end{eqnarray*}
Therefore
\begin{eqnarray*}
l(y, T)-l(x, \tau)&\geq &-{1\over 4}\int_t^T e^{2Kt}{\dot s}^2(t)dt-\int_\tau^T {m\over 2t}e^{2Kt}dt\\
&=&-{Kd^2(x, y)\over 2(e^{-2K\tau}-e^{-2KT})}-\int_\tau^T {m\over 2t}e^{2Kt}dt.
\end{eqnarray*}
Note that $\int_\tau^T {e^{2Kt}\over t}dt\leq \log\left({T\over \tau}\right)+e^{2KT}-e^{2K\tau}$. Thus
\begin{eqnarray*}
\log u(y, T)-\log u(x, \tau)\geq -{Kd^2(x, y)\over 2(e^{-2K\tau}-e^{-2KT})}-{m\over 2}\left[\log\left({T\over \tau}\right)+e^{2KT}-e^{2K\tau}\right].
\end{eqnarray*}
Using ${1\over 1-e^{-x}}\leq {1+x\over x}$, we can derive the desired estimate. \hfill $\square$

\subsection{Proof of Theorem \ref{Thm3}}

Let $\psi(t)={1-e^{-Kt}\over K}$, and  $h=\psi \left[Lu+{|\nabla u|^2\over u}\right]-u[m+4\log(A/u)]$. Then 
$$\psi'+K\psi=1.$$ 
By $(\ref{GB4})$ , under the assumption $Ric_{m, n}(L)\geq -K$ we have
\begin{eqnarray*}
(\partial_t-L){|\nabla u|^2\over u}\leq -{2\over m u}\left|L u-{|\nabla u|^2\over u}\right|^2+2K{|\nabla u|^2\over u},
\end{eqnarray*}
which yields
\begin{eqnarray*}
(\partial_t-L)h\leq -{2\psi\over mu}\left|L u-{|\nabla u|^2\over u}\right|^2+\psi' \left[L u-{|\nabla u|^2\over u}\right]-2{|\nabla u|^2\over u}.
\end{eqnarray*}
By analogue of Hamilton\cite{H1}, we can verify that
\begin{eqnarray*}
{\partial h\over \partial t}\leq Lh \ \ \ {\rm whenever}\ \ h\geq 0.
\end{eqnarray*}
Indeed, we can verify this by examining three cases:\\

$(i)$ If $L u\leq {|\nabla u|^2\over u}$, then $(\partial_t-L)h\leq 0$ since $\psi'\geq 0$.

$(ii)$ If ${|\nabla u|^2\over u}\leq Lu\leq 3{|\nabla u|^2\over u}$, then $(\partial_t-L)h\leq 0$ since $\psi'\leq 1$.

$(iii)$ If $3{|\nabla u|^2\over u}\leq Lu$, then whenever $h\geq 0$, we have
\begin{eqnarray*}
2\left[L u-{|\nabla u|^2\over u}\right]\geq L u+{|\nabla u|^2\over u}={h\over \psi}+{mu+4u\log(A/u)\over \psi} \geq {mu\over \psi},
\end{eqnarray*}
which yields, since $\psi'\leq 1$, we have
\begin{eqnarray*}
(\partial_t-L)h\leq (\psi'-1)\left[Lu-{|\nabla u|^2\over u}\right]-2{|\nabla u|^2\over u}\leq 0.
\end{eqnarray*}

Note that $h\leq 0$ at $t=0$. By the weak maximum principle on complete Riemannian manifolds, see e.g. Theorem 12.10 in \cite{CCLLN}, we conclude that $h\leq 0$ for all $t\in [0, T]$. Thus
\begin{eqnarray*}
{Lu\over u}+{|\nabla u|^2\over u^2}\leq {K\over 1-e^{-Kt}}\left[m+4\log(A/u)\right].
\end{eqnarray*}
This completes the proof of Theorem \ref{Thm3}. \hfill $\square$

\subsection{Proof of Theorem \ref{thm4}}

The lower bound estimate  of $(\ref{htt})$ follows from $(\ref{Hamil})$. It remains to prove the upper bound estimate. Recall the following 

\begin{proposition}\label{prop3} (\cite{Li16}) Suppose that there exist some constants $m\geq n$, $m\in \mathbb{N}$ and $K\geq 0$ such that $Ric_{m, n}(L)\geq -K$. Then, for any small $\varepsilon>0$, there exist some constants $C_i=C_i(m, n, K, \varepsilon)>0$, $i=1, 2$, such that for all $x, y\in M$ and $t>0$,
\begin{eqnarray*}
p_t(x, y)&\leq& {C_1\over \mu(B_y(\sqrt{t}))}\exp\left(-{d^2(x, y)\over 4(1+\varepsilon)t}+\alpha \varepsilon Kt\right)\nonumber\\
& &\hskip1cm \times \left({d(x, y)+\sqrt{t}\over \sqrt{t}}\right)^{m/2} \exp\left({\sqrt{(m-1)K}d(x, y)\over 2}\right),\label{heat-upp}
\end{eqnarray*}
where $\alpha$ is a constant depending only on $m$, and
\begin{eqnarray*}
p_t(x, y)\geq C_2e^{-(1+\varepsilon)\lambda_{K,
m}t}\mu^{-1}(B_y(\sqrt{t}))\exp\left(-{d^2(x, y)\over
4(1-\varepsilon)t}\right)\left[{\sqrt{K}d(x, y)\over \sinh \sqrt{K}d(x,
y)}\right]^{m-1\over 2},\label{heat-lower}
\end{eqnarray*}
where
$$\lambda_{K, m}={(m-1)^2K\over 8}.$$
\end{proposition}

Fix $T>0$,  and let $u(t, x)$ be a positive and bounded solution to the heat equation
$\partial_t u=L u$, $t\in (0, t_1)$. Let
$$
A:=\sup\limits\{u(t, x): 0\leq t\leq t_1, x\in M\}.$$ By
Hamilton's Harnack inequality $(\ref{HHHH4})$, we have
\begin{eqnarray}
t \partial_t \log u(x, t)\leq (1+Kt)\left[m+4\log(A/u(t, x))\right], \ \ \ \forall
(t, x)\in [0, t_1]\times M. \label{EE1}
\end{eqnarray}
Let $s\in (0, T]$, $y\in M$, $t_1=s/2$ and $u(t, x)=p_{s/2+t}(x,
y)$. By $(\ref{EE1})$ and using the upper bound and lower bound estimates of the heat kernel $p_t(x, y)$ in Proposition \ref{prop3},  we have
\begin{eqnarray*}
& &{t\over 2}\partial_t \log p_{s/2+t}(x, y)\leq C_{K, m, T}(1+Kt/2) \\
& &\hskip1cm \times \left[1+{d(x, y)\over \sqrt{t}}+\log\left({C_1\over
C_2}{\mu(B(y, \sqrt{s/2+t}))\over \mu(B(y, \sqrt{s/2}))} \exp\left(C_3{d^2(x,
y)\over s/2+t}+C_4d(x, y)\right)\right)\right].
\end{eqnarray*}
In particular, taking $t=s/2$ and changing $s$ by $t$ we get
\begin{eqnarray*}
& &{t\over 2} t\partial_t \log p_{t}(x, y) \leq  C_{K, m, T}(1+Kt/2)\\
& &\hskip1cm \times \left[1+{d(x, y)\over \sqrt{t}}+\log\left({C_1\over C_2}{\mu(B(y, \sqrt{t}))\over \mu(B(y, \sqrt{t/2}))}
\exp\left(C_3{d^2(x, y)\over t}+C_4d(x, y)\right)\right)\right].
\end{eqnarray*}
By the generalized Bishop-Gromov volume comparison theorem for weighted volume measure, see \cite{Qian, Lot, Li05, WW}, as $Ric_{m, n}(L)\geq -K$, for all $R>r>0$ and $y\in M$, we have 
\begin{eqnarray*}
{\mu(B(y, R))\over \mu(B(y, r))}\leq \left({R\over r}\right)^{m}\exp\left(\sqrt{(m-1)K}R\right).
\end{eqnarray*}
It follows that
\begin{eqnarray*}
t\partial_t \log p_t(x, y)\leq C_{K, m, T}\left(1+{d^2(x,
y)\over t}+{d(x, y)\over \sqrt{t}}+d(x, y)\right),
\end{eqnarray*}
which yields
\begin{eqnarray*}
\partial_x \log p_{t}(x, y)\leq C_{K, m, T}\left(1+{d(x, y)\over
t}+{1\over \sqrt{t}}\right)^2.
\end{eqnarray*}
This completes the proof of Theorem \ref{thm4}. \hfill $\square$

\section{$W$-entropy for Witten Laplacian with $CD(-K, m)$-condition}

Recall that, Perelman \cite{P1} introduced the
$W$-entropy and proved its monotonicity along the conjugate heat
equation associated to the Ricci flow. In \cite{N1, N2}, Ni proved
the monotonicity of the $W$-entropy for the heat equation of the
usual Laplace-Beltrami operator on complete Riemannian manifolds
with non-negative Ricci curvature. In \cite{Li12, Li16},  the second
author of this paper proved the $W$-entropy formula  and its monotonicity and rigidity theorems for the heat
equation of the Witten Laplacian on complete Riemannian manifolds satisfying the $CD(0, m)$-condition
and gave a probabilistic interpretation of the $W$-entropy for the
Ricci flow. In \cite{LL15}, we gave a new proof of the $W$-entropy
formula obtained in \cite{Li12} for the Witten Laplacian  by
using Ni's $W$-entropy formula to the
Laplace-Beltrami operator on $M\times S^{m-n}$ equipped with a suitable warped product Riemannian metric,
and further proved the monotonicity of the $W$-entropy for the heat
equation of the time dependent Witten Laplacian on compact
Riemannian manifolds equipped with the super Ricci flow with respect
to the $m$-dimensional Bakry-Emery Ricci curvature. 
As we have already seen in Section $1$, there is a close connection between the Perelman $W$-entropy for the heat equation of the Witten Laplacian and the Li-Yau Harnack inequality $(\ref{LYm})$ 
on complete Riemannian manifolds satisfying the $CD(0,
m)$-condition.  In this section, we will introduce the Perelman $W$-entropy and prove its monotonicity for the heat equation of the Witten Laplacian on complete Riemannian manifolds with the $CD(-K, m)$-condition. 

Recall the following entropy dissipation formulas  for the heat equation of the Witten Laplacian on complete Riemannian manifolds with bounded geometry condition. In the case of  compact Riemannian manifolds, it is a well-known result due to Bakry and Emery \cite{BE}.

\begin{theorem}\label{AAA}(\cite{Li12, Li16, LL15})
Let~$(M,g)$~be a complete Riemannian manifold with bounded geometry condition, and $\phi\in C^4(M)$ such that $\nabla^k\phi$ are uniformly bounded  on $M$ for $1\leq k\leq 4$. Let $u$ be the fundamental solution to the heat equation ~$\partial_tu=L u$.
Let
$$H(u(t))=-\int_M u\log ud\mu.$$ Then\footnote{The first order entropy dissipation formula $(\ref{ENTH1})$ holds if $\phi\in C^2(M)$ such that $Ric_{m, n}(L)\geq -K$.}
\begin{eqnarray}
{d\over dt}H(u(t))&=&\int_M |\nabla \log  u|^2 ud\mu, \label{ENTH1}\\
{d^2\over dt^2} H(u(t))&=&-2\int_M \Gamma_2(\nabla \log u, \nabla \log u) ud\mu, \label{ENTH2}
\end{eqnarray}
where 
\begin{eqnarray*}
\Gamma_2(\nabla \log u, \nabla \log u)=|\nabla^2\log u|^2+Ric(L)(\nabla \log u, \nabla \log u).
\end{eqnarray*}
\end{theorem}

Let $(M, g, \phi)$ be as in Theorem \ref{AAA}. Inspired by \cite{P1, N1, Li12, Li16, LL15},
we define
\begin{eqnarray*}
H_{m, K}(u, t)=-\int_M u\log u d\mu-\Phi_{m, K}(t),
\end{eqnarray*}
where $\Phi_{m, K}\in C((0, \infty), \mathbb{R})$ satisfies
\begin{eqnarray*}
\Phi_{m, K}'(t)={m\over 2t}e^{4Kt}, \ \ \ \forall t>0.
\end{eqnarray*}

\begin{proposition} Let $(M, g)$ be a complete Riemannian manifold with bounded geometry
condition, $\phi\in C^4(M)$ be such that $\nabla^k\phi$ are uniformly bounded on $M$ for $1\leq k\leq 2$.
Then, under the condition $Ric_{m, n}(L)\geq -K$, we have
\begin{eqnarray*}
{d\over dt}H_{K, m}(u, t)\leq 0.
\end{eqnarray*}
\end{proposition}
{\it Proof}. By the entropy dissipation formulas in Theorem \ref{AAA} and using the fact $\int_M
\partial_t u d\mu=\int_M Lu d\mu=0$, we have
\begin{eqnarray*}
{d\over dt}H_{m, K}(u, t)&=&\int_M \left[{|\nabla u|^2\over u^2}-{m\over 2t}e^{4Kt}\right] ud\mu\label{HmK1}\\
&=&\int_M \left[{|\nabla u|^2\over u^2}-{m\over 2t}e^{4Kt}-e^{2Kt}{\partial_t u\over u}\right] ud\mu.\label{HmK2}
\end{eqnarray*}
By the Hamilton Harnack inequality in Theorem \ref{HLYH}, we have
\begin{eqnarray*}
{d\over dt}H_{m, K}(u, t)\leq 0.
\end{eqnarray*}
\hfill $\square$

\begin{proposition}\label{prop2} Under the same condition as in
Theorem \ref{Th-W2}, we have
\begin{eqnarray*}
{d^2\over dt^2}H_{m, K}(u, t)=-2\int_M [|\nabla^2\log u|^2+Ric(L)(\nabla\log u, \nabla \log u)]ud\mu-\left({2mK\over t}-{m\over 2t^2}\right)e^{4Kt}.
\end{eqnarray*}
\end{proposition}
{\it Proof}. Indeed, by the second order dissipation formula of the Boltzmann entropy in Theorem \ref{AAA}, we have
\begin{eqnarray*}
{d\over dt}\int_M {|\nabla u|^2\over u}d\mu=-2\int_M [|\nabla^2\log u|^2+Ric(L)(\nabla\log u, \nabla \log u)]ud\mu.
\end{eqnarray*}
Combining this with $(\ref{HmK1})$, Proposition \ref{prop2} follows.
\hfill $\square$

\medskip

Based on the Hamilton differential Harnack inequality $(\ref{Hamil})$ in Theorem \ref{HLYH}, we now
 introduce the $W$-entropy for the heat equation $(\ref{HLu})$ of the Witten Laplacian on complete Riemannian manifolds with the $CD(-K, m)$-condition as
follows
\begin{eqnarray*}
W_{m, K}(u, t)={d\over dt}( tH_{m, K}(u, t)).
\end{eqnarray*}
By the entropy dissipation formula in Theorem \ref{AAA}, we have
\begin{eqnarray*}
W_{m, K}(u, t)&=&\int_M \left[t(|\nabla \log u|^2-\Phi'_{m,
K}(t))-\log u-\Phi_{m, K}(t)\right]ud\mu\\
&=&\int_M \left[t(2L(-\log u)-|\nabla \log u|^2)-\log u-\Phi_{m,
K}(t)-\Phi_{m, K}'(t)\right]ud\mu.
\end{eqnarray*}

We are now in a position to prove the main result of this section, i.e., Theorem \ref{Th-W2}.

%
%

\noindent{\bf Proof of Theorem \ref{Th-W2}}. By $(\ref{HmK1})$ and Proposition \ref{prop2}, we have
\begin{eqnarray*}
{d\over dt}W_{m, K}(u, t)
&=&-2t\left[\int_M [|\nabla^2\log u|^2+Ric(L)(\nabla\log u, \nabla \log u)]ud\mu+\left({mK\over t}-{m\over 4t^2}\right)e^{4Kt}\right]\\
& &\hskip2cm \ \ +2\int_M \left[{|\nabla u|^2\over u^2}-{m\over 2t}e^{4Kt}\right] ud\mu.
\end{eqnarray*}
Note that
\begin{eqnarray*}
\left|\nabla^2\log u+\left({e^{2Kt}\over
2t}+a(t)\right)g\right|^2=|\nabla^2\log u|^2+2\left({e^{2Kt}\over
2t}+a(t)\right)\Delta \log u+n\left({e^{2Kt}\over 2t}+a(t)\right)^2.
\end{eqnarray*}
By direct calculation, we have
\begin{eqnarray*}
{d\over dt}W_{m, K}(u, t)
&=&-2t\int_M \left|\nabla^2\log u+\left({e^{2Kt}\over 2t}+a(t)\right)g\right|^2 u d\mu\\
& &-2t\int_M \left(Ric_{m, n}(L)+\left(2a(t)-{1-e^{2Kt}\over t}\right)g\right)(\nabla\log u, \nabla \log u) ud\mu\\
& &\ \ \ +2nt\left({e^{2Kt}\over 2t}+a(t)\right)^2-{me^{4Kt}\over 2t}-2mK e^{4Kt}\\
& &\ \ \ +2(e^{2Kt}+2ta(t))\int_M \nabla \phi\cdot \nabla\log  u \
ud\mu-2t\int_M {|\nabla \phi\cdot\nabla \log u|^2 \over m-n}ud\mu.
\end{eqnarray*}
Let $a(t)$ be chosen such that $2a(t)-{1-e^{2Kt}\over t}=K$.
Then
\begin{eqnarray*}
{d\over dt}W_{m, K}(u, t)
&=&-2t\int_M \left[\left|\nabla^2\log u+\left({K\over 2}+{1\over 2t}\right)g\right|^2+(Ric_{m, n}(L)+Kg)(\nabla\log u, \nabla \log u)\right] ud\mu\\
& &\ \ \ +2nt\left({1\over 2t}+{K\over 2}\right)^2-{me^{4Kt}\over 2t}-2mK e^{4Kt}\\
& &\ \ \ +2(1+Kt)\int_M \nabla \phi\cdot \nabla\log  u \ ud\mu-2t\int_M {|\nabla \phi\cdot\nabla \log u|^2 \over m-n}ud\mu.
\end{eqnarray*}
Combining this with
\begin{eqnarray*}
& &{1\over m-n}\int_M \left|\nabla\phi\cdot\nabla\log u-{(m-n)(1+Kt)\over 2t}\right|^2ud\mu\\
&=&{(m-n)(1+Kt)^2\over 4t^2}-{1+Kt\over t}\int_M \nabla \phi\cdot \nabla\log  u \ ud\mu+\int_M {|\nabla \phi\cdot\nabla \log u|^2 \over m-n}ud\mu,
\end{eqnarray*}
and noting that
\begin{eqnarray*}
& & 2nt\left({1\over 2t}+{K\over 2}\right)^2-{me^{4Kt}\over 2t}-2mK e^{4Kt}+{(m-n)(1+Kt)^2\over 2t}\\
& &\hskip3cm ={m\over 2t}\left[(1+Kt)^2-e^{4Kt}(1+4Kt)\right],
\end{eqnarray*}
we can derive the desired  $W$-entropy formula. The rest of the
proof is obvious.  \hfill $\square$

In particular, taking $m=n$, $\phi\equiv 0$ and $g$ is a fixed Riemannian metric, we have the following $W$-entropy formula for the heat equation of the Laplace-Beltrami operator on Riemannian manifolds, which extends Ni's result in \cite{N1} for $K=0$.

\begin{theorem} \label{Th-W2c} Let $(M, g)$ be a complete Riemannian manifold with bounded geometry
condition. Let $u$ be the fundamental solution to the heat equation
$\partial_t u=\Delta u$. Then
\begin{eqnarray*}
{d\over dt}W_{n, K}(u, t)
&=&-2t\int_M \left[\left|\nabla^2\log u+\left({K\over 2}+{1\over 2t}\right)g\right|^2+(Ric+Kg)(\nabla\log u, \nabla \log u)\right] ud\mu\\
& &\hskip3cm -{n\over 2t}\left[e^{4Kt}(1+4Kt)-(1+Kt)^2\right].
\end{eqnarray*}
In particular, if $Ric\geq -K$, then, for all $t\geq 0$, we have
\begin{eqnarray*}
{d\over dt}W_{n, K}(u, t)\leq -{n\over
2t}\left[e^{4Kt}(1+4Kt)-(1+Kt)^2\right].
\end{eqnarray*}
Moreover, the equality holds at some time $t=t_0>0$ if and only if
$M$ is an Einstein manifold, i.e., $Ric=-Kg$, and the potential
function $f=-\log u$ satisfies the shrinking soliton equation, i.e.,
\begin{eqnarray*}
Ric+2\nabla^2f={g\over t}.
\end{eqnarray*}
\end{theorem}

By analogue of the $W$-entropy for the heat equation of the Witten Laplacian on complete Riemannian manifolds with the $CD(-K, m)$-condition, we can prove the $W$-entropy formula for the heat equation of the time dependent 
Witten Laplacian on compact manifolds equipped with a $(-K, m)$-super Ricci flow. To do so, let us recall the entropy dissipation formula on compact manifolds with time dependent metrics and potentials.

\begin{theorem}\label{LLPJM} (\cite{LL15})  Let $(M, g(t), t\in [0, T])$ be a family of compact Riemannian manifolds with potential
functions
$\phi(t)\in C^\infty(M)$, $t\in [0, T]$. Suppose that $g(t)$ and $\phi(t)$
satisfy the conjugate
equation
\begin{eqnarray*}
\frac{\partial \phi} {\partial t}={1\over 2}{\rm Tr}\left(
\frac{\partial g}{\partial t}\right).
\end{eqnarray*}
Let
\begin{eqnarray*}
L=\Delta_{g(t)}-\nabla_{g(t)}\phi(t)\cdot\nabla_{g(t)}
\end{eqnarray*}
be the time dependent Witten Laplacian on $(M, g(t), \phi(t))$. Let $u$ be a positive solution of the heat equation
\begin{eqnarray*}
\partial_t u = Lu
\end{eqnarray*}
with initial data $u(0)$ satisfying $\int_M
u(0)d\mu(0)=1$.
Let
$$H(u, t)=-\int_M u\log u d\mu$$
be the Boltzmann-Shannon entropy for the heat equation $\partial_t u=Lu$. Then
\begin{eqnarray*}
{\partial \over \partial t } H(u, t)&=& \int_M |\nabla \log u|^2_{g(t)} ud\mu,\label{1stH} \\
{\partial^2\over \partial t^2} H(u, t)&=&-2\int_M \left[|\nabla^2\log u|^2+\left({1\over 2}{\partial g\over \partial t}+Ric(L)\right)(\nabla \log u, \nabla \log u)\right]u d\mu.\label{2ndH}
\end{eqnarray*}
\end{theorem}

\noindent{\bf Proof of Theorem \ref{Th-W3}}. Base on the entropy dissipation formulas in Theorem \ref{LLPJM}, the proof  of Theorem \ref{Th-W3} is similar to the one of Theorem \ref{Th-W2}. See \cite{LL15} for the case $K=0$. \hfill $\square$

\section{Comparison with another $W$-entropy functional}

To end this paper, let us mention that in our previous paper \cite{LL15} we introduced another $W$-entropy functional for  the heat equation associated with the Witten Laplacian on complete Riemannian manifolds satisfying the $CD(-K, m)$-condition as follows

\begin{eqnarray*}
\widetilde{W}_{m, K}(u)={d\over dt}(t\widetilde{H}_{m, K}(u)),
\end{eqnarray*}
where
\begin{eqnarray*}
\widetilde{H}_{m, K}(u)=-\int_M u\log ud\mu-\left[{m\over 2t}(1+\log(4\pi t))+{mKt\over 2}(1+{1\over 6}Kt)\right],
\end{eqnarray*}
and we proved that
\begin{eqnarray*}
{d\over dt}\widetilde{W}_{m, K}(u)&=&-2t\int_M \left[\left|\nabla^2\log u+\left({K\over 2}+{1\over 2t}\right)g\right|^2+(Ric_{m, n}(L)+Kg)(\nabla\log u, \nabla \log u)\right] ud\mu\\
& &\hskip1.5cm -{2t\over m-n}\int_M \left|\nabla \phi\cdot \nabla\log  u-{(m-n)(1+Kt)\over 2t}\right|^2ud\mu.
\end{eqnarray*}

It is interesting to compare the $W$-entropy defined in \cite{LL15} with the $W$-entropy defined in this paper,  and to compare the $W$-entropy formula proved in \cite{LL15} with the $W$-entropy formula obtained  in Theorem \ref{Th-W2}. Indeed, 
letting 
\begin{eqnarray*}
\Psi_{m, K}(t)=\Phi_{m, K}(t)-\left[{m\over 2t}(1+\log(4\pi t))+{mKt\over 2}(1+{1\over 6}Kt)\right],\label{Psi1}
\end{eqnarray*} we have
\begin{eqnarray*}
\widetilde{W}_{m, K}(u)-W_{m, K}(u)={d\over dt}(t\Psi_{m, K}(t)), \label{Psi2}
\end{eqnarray*}
Moreover,  by direct calculation we have 
\begin{eqnarray*}
{d\over dt}(\widetilde{W}_{m, K}(u)-W_{m, K}(u))&=&{d^2\over dt^2}(t\Psi_{m, K}(t))\nonumber\\
&=&{m\over 2t}\left[e^{4Kt}(1+4Kt)-(1+Kt)^2\right].\label{Psi3}
\end{eqnarray*}
This explains clearly the difference between the $W$-entropy defined in \cite{LL15}  and the $W$-entropy defined in this paper, and the difference between the $W$-entropy formula proved in \cite{LL15} with the $W$-entropy formula obtained  in Theorem \ref{Th-W2}. 

Similarly, we can reformulate Theorem \ref{Th-W3} in terms of $\widetilde {W}_{m, K}$. See \cite{LL15, LL17c}.

\medskip

\noindent{\bf Acknowledgement}.  Part of this work was done when the second author visited the Institut des Hautes Etudes Scientifiques and the Max-Planck Institute for Mathematics Bonn. The authors would like to thank Professors D. Bakry, J.-M. Bismut, M. Ledoux, N. Mok, K.-T. Sturm, A. Thalmaier, F.-Y. Wang and Dr. Yuzhao Wang for their 
interests and helpful discussions during  the preparation of this paper.  We are very grateful to anonymous referee for his careful reading and for his very nice comments which lead us to improve the writting of this paper.

\medskip

\begin{flushleft}
\medskip\noindent

Songzi Li, School of Mathematical Science, Beijing Normal University, No. 19, Xin Jie Kou Wai Da Jie, 100875, China, Email: songzi.li@bnu.edu.cn

\medskip

Xiang-Dong Li, Academy of Mathematics and Systems Science, Chinese
Academy of Sciences, 55, Zhongguancun East Road, Beijing, 100190,  China, 
E-mail: xdli@amt.ac.cn

and

School of Mathematical Sciences, University of Chinese Academy of Sciences, Beijing, 100049, China
\end{flushleft}


\begin{thebibliography}{99}
\bibitem{BE} D. Bakry, M. Emery, Diffusion hypercontractives, S\'em. Prob. XIX, Lect. Notes in Maths. 1123 (1985), 177-206.
\bibitem{BL} D. Bakry, M. Ledoux, A logarithmic Sobolev form of the Li-Yau parabolic inequality, Rev. Mat. Iberoam. 22 (2006), No. 2, 683-702.
\bibitem{Cal58} E. Calabi, An extension of E. Hopf's maximum principle with an application to Riemannian geometry, Duke Math. J. 25 (1958), 1, 45-56.
\bibitem{CLN} B. Chow, P. Lu, L. Ni, {\it Hamilton's Ricci Flow, Lectures
in Comtemporary Maths.} Sciences Press, Beijing, Amer. Math. Soc.
2006.
\bibitem{CCLLN} B. Chow, S.-C. Chu,  D. Glickenstein,  C. Guenther,  J. Isenberg,  T. Ivey, D. Knopf,  P. Lu,  F. Luo,  L. Ni, {\it  The Ricci Flow: Techniques and Applications: Part II: Analytic Aspects}, Mathematical Surveys and Monographs
Volume 144,  2008, 458 pp. 
\bibitem{H1} R.S. Hamilton, A matrix Harnack estimate for the heat equation, Comm. Anal.
Geom. 1 (1993), no. 1, 113-126.
\bibitem{St4} E. Kopfer, K.-T. Sturm, Heat flows on time-dependent metric measure spaces and super-Ricci flows, arXiv:1611.02570 
\bibitem{LY} P. Li, S.-T. Yau, On the parabolic kernel of the
Schr\"odinger operator, Acta. Math. 156 (1986), 153-201.
\bibitem{Li05}X.-D. Li, Liouville theorems for symmetric diffusion
operators on complete Riemannian manifolds, J. Math. Pures Appl. 84 (2005), 1295-1361.
\bibitem{Li12}X.-D. Li, Perelman's entropy formula for the Witten Laplacian
on Riemannian manifolds via Bakry-Emery Ricci
curvature, Math. Ann. 353 (2012), 403-437.
\bibitem{Li11}X.-D. Li, Perelman's $W$-entropy for the Fokker-Planck equation
over complete Riemannian manifolds, Bull. Sci. Math. 135 (2011) 871-882.
\bibitem{Li16} X.-D. Li, Hamiltons Harnack inequality and the W-entropy
formula on complete Riemannian manifolds, Stoch. Processes and Appl., 126 (2016) 1264-1283.
\bibitem{LL15}S. Li, X.-D. Li,  $W$-entropy formula for the Witten Laplacian on
manifolds with time dependent metrics and potentials,  Pacific J. Math. Vol. 278 (2015), No. 1, 173-199.
\bibitem{LL14} S. Li, X.-D. Li, Harnack inequalities and $W$-entropy formula for Witten Laplacian on Riemannian manifolds with $K$-super Perelman Ricci flow, arxiv1412.7034,  version1 (22 December 2014) and version 2 (7 February 2016). 
\bibitem{LL17a} S. Li, X.-D. Li, On Harnack inequalities for Witten Laplacian on Riemannian manifolds with super Ricci flows, to appear in Special Issue of Asian J. Math. in honor of Prof. N. Mok's 60th birthday, arXiv:1706.05304
\bibitem{LL17b} S. Li, X.-D. Li, $W$-entropy, super Perelman Ricci flows and $(K, m)$-Ricci solitons,  arXiv:1706.07040
\bibitem{LL17c} S. Li, X.-D. Li, $W$-entropy formulas on super Ricci flows and Langevin deformation on Wasserstein space over Riemannian manifolds, submitted to Science China Mathematics, 2017. 
\bibitem{LX} J. Li, X. Xu, Differential Harnack inequalities on Riemannian manifolds I: linear heat equation, Adv. Math. 226 (5) (2011), 4456-4491.
\bibitem{Lot} J. Lott, Some geometric properties of the Bakry-Emery Ricci tensor. Comment. Math. Helv. {\bf 78} (2003), 865-883.
\bibitem{N1} L. Ni, The entropy formula for linear equation, J. Geom. Anal. 14 (1), 87-100, (2004).
\bibitem{N2} L. Ni, Addenda to ``The entropy formula for linear equation'', J. Geom. Anal. 14 (2), 329-334, (2004).
\bibitem{P1} G. Perelman, The entropy formula for the Ricci flow and its geometric applications, http://arXiv.org/abs/maths0211159.
\bibitem{Qian} Z. Qian, Estimates for weighted volumes and applications. Quart. J. Math. Oxford Ser. (2), 48(190), 235-242, 1997.
\bibitem{St3} K.-T. Sturm,  Super-Ricci flows for metric measure spaces, arXiv:1603.02193
\bibitem{WW} G. Wei, W. Wylie, Comparison geometry for the Bakry-Emery Ricci tensor, J. Diff. Geom. (2009), 377-405.
\end{thebibliography}
\end{document}